\documentclass[10pt,twocolumn]{article}

\usepackage{amssymb,euscript,stmaryrd,epsfig}

\usepackage{amssymb,amsmath}
\usepackage{epsfig}
\usepackage{amsfonts}
\usepackage{latexsym}
\usepackage{xspace}
\usepackage{graphicx}
\usepackage{color}

\usepackage{subfigure}

 \usepackage{setspace}

\newcommand{\R}{\mathbb R}

\newcommand{\N}{\mathbb N}

\title{Domain decomposition methods for compressed sensing}

\author{Massimo Fornasier\thanks{Johann Radon Institute for Computational and Applied Mathematics (RICAM),
Austrian Academy of Sciences, Altenbergerstrasse 69, A-4040, Linz, Austria Email: {\tt massimo.fornasier@oeaw.ac.at}} \and Andreas Langer\thanks{Johann Radon Institute for Computational and Applied Mathematics (RICAM),
Austrian Academy of Sciences, Altenbergerstrasse 69, A-4040, Linz, Austria Email: {\tt andreas.langer@oeaw.ac.at}} \and Carola-Bibiane Sch\"onlieb\thanks{Department of Applied Mathematics and Theoretical Physics (DAMTP),
Centre for Mathematical Sciences,
Wilberforce Road,
Cambridge CB3 0WA,
United Kingdom.Email: {\tt c.b.s.schonlieb@damtp.cam.ac.uk} } }

\begin{document}

\graphicspath{{./graphics/}}

\maketitle

\pagestyle{myheadings}
\thispagestyle{plain}
\markboth{M. FORNASIER, A. LANGER, AND C.-B. SCH\"ONLIEB}{Domain decomposition methods for compressed sensing}

\begin{abstract}
We present several domain decomposition algorithms for sequential and parallel minimization of functionals formed by a discrepancy term with respect to data and total variation constraints. The convergence properties of the algorithms are analyzed. We provide several numerical experiments, showing the successful application of the algorithms for the restoration 1D and 2D signals in interpolation/inpainting problems respectively, and in a compressed sensing problem, for recovering piecewise constant medical-type images from partial Fourier ensembles.
\end{abstract}

\section{Introduction}

In concrete applications for image processing, one might be interested to recover at best a digital image provided only partial linear or nonlinear measurements, possibly corrupted by noise.
Given the observation that natural and man-made images are characterized by a relatively small number of edges and extensive relatively uniform parts, one may want to help the reconstruction by imposing that the interesting solution is the one which matches the given data and has also a few discontinuities localized on sets of lower dimension.

In the context of compressed sensing \cite{carotaXX,do04}, it has been clarified the fundamental role of minimizing $\ell_1$-norms in order to promote sparse solutions. This understanding furnishes an important interpretation of total variation minimization, i.e., the minimization of the $\ell_1$-norm of derivatives \cite{ROF}, as a regularization technique for image restoration. 
Several numerical strategies to perform efficiently total variation minimization have been proposed in the literature. We list a few of the relevant ones, ordered by their chronological appearance:

(i) the approach of Chambolle and Lions \cite{ChL} by re-weighted least squares, see also \cite{DDFG} for generalizations and refinements in the context of compressed sensing;

(ii) variational approximation via local quadratic functionals as in the work of Vese et al. \cite{Ve01,AK02};

(iii) iterative thresholding algorithms based on projections onto convex sets as in the work of Chambolle \cite{Ch} as well as in the work of Combettes-Wajs \cite{CW} and Daubechies et al. \cite{dateve06};

(iv)  iterative minimization of the Bregman distance as in the work of Osher et al. \cite{breg};

(v)  the approach proposed by Nesterov \cite{ne05} and its modifications by Weiss et al. \cite{WBA}.

These approaches differ significantly, and it seems that the ones collected in the groups iv) and v) do show presently the best performances in practice.
However, none of the mentioned methods is able to address in real-time, or at least in an acceptable computational time, extremely large problems, such as 4D imaging (spatial plus temporal dimensions) from functional magnetic-resonance in nuclear medical imaging, astronomical imaging or global terrestrial seismic tomography.
For such large scale simulations we need to address methods which allow us to reduce the problem to a finite sequence of sub-problems of more manageable size, perhaps by one of the methods listed above. With this aim we introduced subspace correction and domain decomposition methods  both for $\ell_1$-norm and total variation minimizations \cite{fo07,fosc08}.
Due to the nonadditivity of the total variation with respect to a domain decomposition (the total variation of a function on the whole domain equals the sum of the total variations on the subdomains plus the size of the jumps at the interfaces), one encounters additional difficulties in showing convergence of such decomposition strategies to global minimizers.

In this paper we review concisely both nonoverlapping and overlapping domain decomposition methods for total variation minimization and we provide their properties of convergence to global minimizers. Moreover, we show numerical applications in classical problems of signal and image processing, such as signal interpolation and image inpainting. We further include applications in the context of compressed sensing for recovering piecewise constant medical-type images from partial Fourier ensembles \cite{carotaXX}. 

\section{Notations and preliminaries} Since we are interested to a discrete setting, we define the domain of our multivariate digital signal $\Omega = \{x_1^1 < \ldots < x_{N_1}^1\} \times \ldots \times \{x_1^d < \ldots < x_{N_d}^d\}\subset\mathbb{R}^d$, $d\in\mathbb{N}$ and we consider the signal space  $\mathcal H = \R^{N_1\times N_2 \times \ldots \times N_d}$, where $N_i\in\N$ for $i=1,\ldots,d$. For $u\in\mathcal{H}$ we write $u=u(x_{\mathbf{i}})_{\mathbf{i}\in\mathcal{I}}$ with index set $\mathcal{I}:=\prod_{k=1}^d\{1,\ldots,N_k\}$ and  $u(x_{\mathbf{i}}) = u(x_{i_1}^1, \ldots , x_{i_d}^d)$ where $i_k\in\{1,\ldots,N_k\}$ and $x_{\mathbf{i}}  \in \Omega$. Then, for $g \in\mathbb R^N$, the $\ell_p$-norm is given by
$
\|g \|_p=\left(\sum_{k=1}^N|g_k|^p\right)^{1/p}, \quad 1 \leq p < \infty$, 
and for $u \in \mathcal H$ the discrete gradient $\nabla u$ is given by $(\nabla u)(x_{\mathbf{i}}) = ((\nabla u)^1(x_{\mathbf{i}}),\ldots,(\nabla u)^d(x_{\mathbf{i}}))$
with  $(\nabla u)^j(x_{\mathbf{i}})=u(x^1_{i_1},\ldots,x^j_{i_j+1},\ldots,x^d_{i_d}) - u(x^1_{i_1},\ldots,x^j_{i_j},\ldots,x^d_{i_d})$ if $i_j<N_j$, and  $(\nabla u)^j(x_{\mathbf{i}})=0$ if $i_j=N_j$, for all $j=1,\ldots,d$ and for all $\mathbf{i}=(i_1,\ldots,i_d)\in\mathcal{I}$. The total variation of $u\in\mathcal{H}$ in the discrete setting is then defined as
$\left|D(u)\right|(\Omega) = \sum\limits_{\mathbf{i}\in\mathcal{I}} \left|(\nabla u)(x_{\mathbf{i}})\right|$, with $|y|=\sqrt{y_1^2+\ldots+y_d^2}$ for every $y=(y_1,\ldots,y_d)\in\R^d$. We define the scalar product of $u,v \in \mathcal{H}$ as usual, $\langle u,v \rangle_{\mathcal{H}}= \sum_{\mathbf{i}\in\mathcal{I}}u(x_{\mathbf{i}})v(x_{\mathbf{i}})$, and the scalar product of $p,q \in \mathcal{H}^d$ as $
\langle p,q \rangle_{\mathcal{H}^d}= \sum_{\mathbf{i}\in\mathcal{I}}p^1(x_{\mathbf{i}})q^1(x_{\mathbf{i}}) + \ldots + p^d(x_{\mathbf{i}})q^d(x_{\mathbf{i}})$.
Further we introduce a discrete divergence $\operatorname{div}: \mathcal{H}^d \to \mathcal{H}$ defined, by analogy with the continuous setting, by $\operatorname{div} = -\nabla^*$ ($\nabla^*$ is the adjoint of the gradient $\nabla$). 
In the following we denote with $\pi_K$ the orthogonal projection onto a closed convex set $K$.

\subsection{Projections onto convex sets}
With these notation, we define the closed convex set 
\begin{eqnarray*}
K:=\left\{\operatorname{div} p: p\in \mathcal H^d, \left|p_{\mathbf i}\right|\leq 1\;\mbox{ for all } \mathbf i \in \mathcal I \right\}.
\end{eqnarray*}

We briefly recall here  an algorithm proposed by Chambolle in \cite{Ch} in order to compute the projection onto $\alpha K$. The following semi-implicit gradient descent algorithm is given to approximate $\pi_{\alpha K}(g)$: Choose $\tau>0$, let $p^{(0)}=0$ and, for any $n\geq 0$, iterate
\begin{eqnarray}\label{chprojit}
p_{\mathbf i}^{(n+1)} = \frac{p_{\mathbf i}^{(n)} + \tau (\nabla(\operatorname{div} p^{(n)}-g/\alpha))_{\mathbf i}}{1+\tau \left|(\nabla(\operatorname{div} p^{(n)}-g/\alpha))_{\mathbf i}\right|}.
\end{eqnarray}
For $\tau>0$ sufficiently small, the iteration $\alpha\operatorname{div} p^{(n)}$ converges to $\pi_{\alpha K}(g)$ as $n\rightarrow\infty$ (compare \cite[Theorem 3.1]{Ch}).

\subsection{Setting of the problem}
Given a model linear operator $T:\mathcal H \to \mathbb R^K$, we are considering the following discrete minimization problem
\begin{equation}\label{problem}
\arg \min_ {u\in\mathcal{H}} \left \{\mathcal{J}(u):=\|Tu-g\|_2^2 + 2\alpha |Du|(\Omega) \right \}
\end{equation}
where $g\in\mathbb R^K$ is a given datum and $\alpha > 0$ is a fixed regularization parameter. Note that, up to rescaling the parameter $\alpha$ and the datum $g$, we can always assume $\|T\| < 1$. Moreover, in order to ensure existence of solutions, we assume $1 \notin \ker(T)$. For both nonoverlapping and overlapping domain decompositions, we will consider a linear sum $\mathcal H = V_1 + V_2$ with respect to two subspaces $V_1$, $V_2$ defined by a suitable decomposition of the physical domain $\Omega$. We restrict our discussion to two subspaces, but the analysis can be extended in a straightforward way to multiple subspaces. With this splitting we want to minimize $\mathcal J$ by suitable instances of the following alternating algorithm:  Pick an initial $V_1 + V_2 \ni  u_1^{(0)}+ u_2^{(0)} : = u^{(0)} \in \mathcal H$, for example $u^{(0)}=0$, and iterate
$$
\left \{ 
\begin{array}{ll}
u_1^{(n+1)} \approx \arg \min_{v_1 \in V_1}  \mathcal  J(v_1 +u_2^{(n)}) &\\
u_2^{(n+1)} \approx  \arg \min_{v_2 \in V_2} \mathcal J(u_1^{(n+1)} + v_2) &\\
u^{(n+1)}:=u_1^{(n+1)} + u_2^{(n+1)}.
\end{array}
\right.
$$

\section{Nonoverlapping domain decomposition methods}

Let us consider the disjoint domain decomposition $\Omega = \Omega_1 \cup \Omega_2$ and $\Omega_1 \cap \Omega_2 = \emptyset$ and the corresponding spaces $V_j=\{ u \in \mathcal H: \operatorname{supp}(u) \subset \Omega_j \}$, for $j=1,2$. Note that $\mathcal H = V_1 \oplus V_2$. It is useful to us to introduce an auxiliary functional $\mathcal J^s_1$, called the {\it surrogate functional} of $\mathcal J$: For $j \in \{1,2\}$ and $\check j \in \{1,2\} \setminus \{j\}$, assume $a,u_j \in V_j$ and $u_{\check j}\in V_{\check j}$ and define
\begin{equation}
\label{surrfunc}
\mathcal J^s_j( u_j+ u_{\check j}, a) := \mathcal J(u_j+ u_{\check j})+ \| u_j -a\|_{\mathcal H}^2 - \| T(u_j -a)\|_{\mathcal H}^2.
\end{equation}
As it will be clarified later, the minimization of $\mathcal J^s_j( u_j+ u_{\check j}, a)$ with respect to $u_j$ and for fixed $u_{\check j}, a$ is an operation which can be realized more easily than the direct minimization of the parent functional $\mathcal J(u_j+ u_{\check j})$ for the sole $u_{\check j}$ fixed.

\subsection{Sequential algorithm}
In the following we denote $u_j = \pi_{V_j} u$ the orthogonal projection onto $V_j$, for $j=1,2$. Let us explicitely express the algorithm  as follows:
 Pick an initial $V_1\oplus V_2 \ni u_1^{(0,L)}+ u_2^{(0,M)} : = u^{(0)} \in \mathcal H$, for example $u^{(0)}=0$, and iterate
\begin{equation}
\label{schw_sp:it2}
\left \{ 
\begin{array}{ll}
\left \{ 
\begin{array}{ll}
u_1^{(n+1,0)} =  u_1^{(n,L)} \mbox{ and for }\ell=0,\dots, L-1 & \\
u_1^{(n+1,\ell+1)} =  \arg \min_{u_1 \in V_1} \mathcal J_1^s(u_1+ u_2^{(n,M)}, u_1^{(n+1,\ell)}) &  \\
\end{array}\right. &\\
\left \{ 
\begin{array}{ll}
u_2^{(n+1,0)} =  u_2^{(n,M)}\mbox{ and for }m=0,\dots, M-1  &\\
u_2^{(n+1,m+1)} = \arg \min_{u_2 \in V_2} \mathcal J_2^s(u_1^{(n+1,L)}+ u_2, u_2^{(n+1,m)}) &\\
\end{array}\right. &\\
u^{(n+1)}:=u_{1}^{(n+1,L)} + u_{2}^{(n+1,M)}.
\end{array}
\right.
\end{equation} 
Note that we do prescribe a finite number $L$ and $M$ of inner iterations for each subspace respectively.

\subsection{Parallel algorithm}

 The parallel version of the previous algorithm reads as follows: 
Pick an initial $V_1\oplus V_2 \ni u_1^{(0,L)}+ u_2^{(0,M)} : = u^{(0)} \in \mathcal H$, for example $u^{(0)}=0$, and iterate
\begin{equation}
\label{schw_sp:it4}
\left \{ 
\begin{array}{ll}
\left \{ 
\begin{array}{ll}
u_1^{(n+1,0)} =  u_1^{(n,L)}\mbox{ and for }\ell=0,\dots, L-1  &\\
u_1^{(n+1,\ell+1)} =  \arg \min_{u_1 \in V_1} \mathcal J_1^s(u_1+ u_2^{(n,M)}, u_1^{(n+1,\ell)}) & \\
\end{array}\right. &\\
\left \{ 
\begin{array}{ll}
u_2^{(n+1,0)} =  u_2^{(n,M)} \mbox{ and for } m=0,\dots, M-1 &\\
u_2^{(n+1,m+1)} = \arg \min_{u_2 \in V_2} \mathcal J_2^s(u_1^{(n,L)}+ u_2, u_2^{(n+1,m)}) &\\
\end{array}\right. &\\
u^{(n+1)}:=\frac{u_{1}^{(n+1,L)} + u_{2}^{(n+1,M)}+u^{(n)}}{2}.
\end{array}
\right.
\end{equation} 
Note that $u^{(n+1)}$ is the average of the current iteration and the previous one as it is the case for successive overrelaxation methods (SOR) in classical numerical linear algebra.

\section{Overlapping domain decomposition methods}

\begin{figure}[htbp]

\begin{picture}(-140,10)(-140,10)

\put(35,23){$\Omega_2$}
\put(-22,+7){\textcolor{red}{$\Gamma_2$}}
\put(-20,20){\line(1,0){120}}
\put(-20,22){\line(0,-1){4}}
\put(100,22){\line(0,-1){4}}

\put(-20,0){\textcolor{red}{\circle{5}}}
\put(20,0){\textcolor{red}{\circle{5}}}
\put(-100,0){\textcolor{blue}{\thicklines\line(1,0){200}}}
\put(-100,-2){\textcolor{blue}{\thicklines\line(0,1){4}}}
\put(100,-2){\textcolor{blue}{\thicklines\line(0,1){4}}}

\put(18,-15){\textcolor{red}{$\Gamma_1$}}
\put(-100,-20){\line(1,0){120}}
\put(20,-22){\line(0,1){4}}
\put(-100,-22){\line(0,1){4}}
\put(-45,-29){$\Omega_1$}
\end{picture}
$\vspace{1cm}$
\caption{\small Overlapping domain decomposition and internal boundaries.}
\label{fig:overlap}
\end{figure}
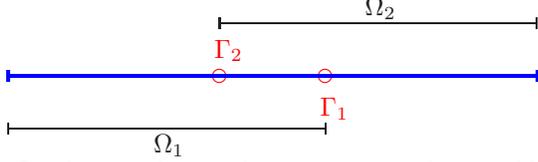

Let us consider the overlapping domain decomposition $\Omega = \Omega_1 \cup \Omega_2$ and $\Omega_1 \cap \Omega_2 \neq \emptyset$ and the corresponding spaces $V_j=\{ u \in \mathcal H: \operatorname{supp}(u) \subset \Omega_j \}$, for $j=1,2$. Note that now $\mathcal H = V_1 + V_2$ is not anymore a direct sum of  $V_1$ and $V_2$, but just a linear sum of subspaces. We define the internal boundaries $\Gamma_j = \partial \Omega_ j \cap \Omega_{\check j}$, $j \in \{1,2\}$ and $\check j \in \{1,2\} \setminus \{j\}$ (see Figure \ref{fig:overlap}).

\subsection{Sequential algorithm}
Associated to the decomposition $\{\Omega_j: j=1,2\}$ let us fix a partition of unity $\{\chi_j:j=1,2\}$, i.e., $\chi_1 + \chi_2 =1$, $\operatorname{supp}(\chi_j) \subset \Omega_j$, and $\chi_j|_{\Gamma_j}=0$.
Let us explicitely express the algorithm now as follows:
Pick an initial $V_1 + V_2 \ni \tilde u_1^{(0)}+ \tilde u_2^{(0)} : = u^{(0)} \in \mathcal{H}$, for example $u^{(0)}=0$, and iterate
\begin{equation}
\label{schw_sp:it21}
\left \{ 
\begin{array}{ll}
\left \{ 
\begin{array}{ll}
u_1^{(n+1,0)} = \tilde u_1^{(n)} \mbox{ and for }\ell=0,\dots, L-1 &\\
u_1^{(n+1,\ell+1)} =  \arg \min_{\stackrel{u_1 \in V_1}{u_1|_{\Gamma_1}=0}} \mathcal J_1^s(u_1+ \tilde u_2^{(n)}, u_1^{(n+1,\ell)}) & \\
\end{array}\right. &\\
\left \{ 
\begin{array}{ll}
u_2^{(n+1,0)} = \tilde{u}_2^{(n)}\mbox{ and for }m=0,\dots, M-1 &\\
u_2^{(n+1,m+1)} = \arg \min_{\stackrel{u_2 \in V_2}{u_2|_{\Gamma_2}=0}} \mathcal J_2^s(u_1^{(n+1,L)}+ u_2, u_2^{(n+1,m)}) &\\
\end{array}\right. &\\
u^{(n+1)}:=u_{1}^{(n+1,L)} + u_{2}^{(n+1,M)}\\
\tilde{u}_1^{(n+1)}:=\chi_1\cdot u^{(n+1)} \\
\tilde{u}_2^{(n+1)}:=\chi_2\cdot u^{(n+1)}
\end{array}
\right.
\end{equation} 

\begin{figure}[ht!]
\begin{center}
\graphicspath{{./}}
    \subfigure[]{\label{l1}\includegraphics[height=3.4cm]{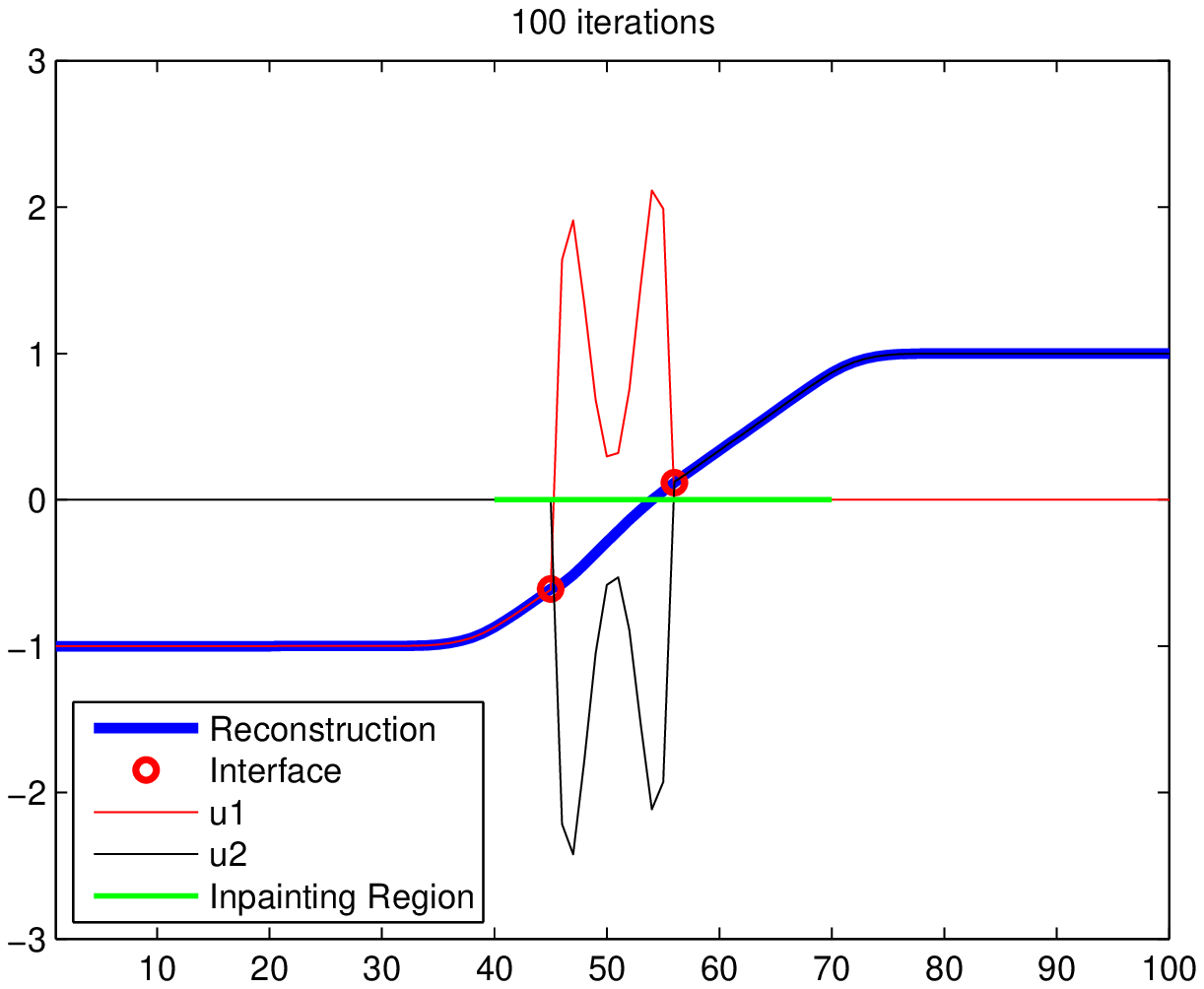}}
    \subfigure[]{\label{l2}\includegraphics[height=3.4cm]{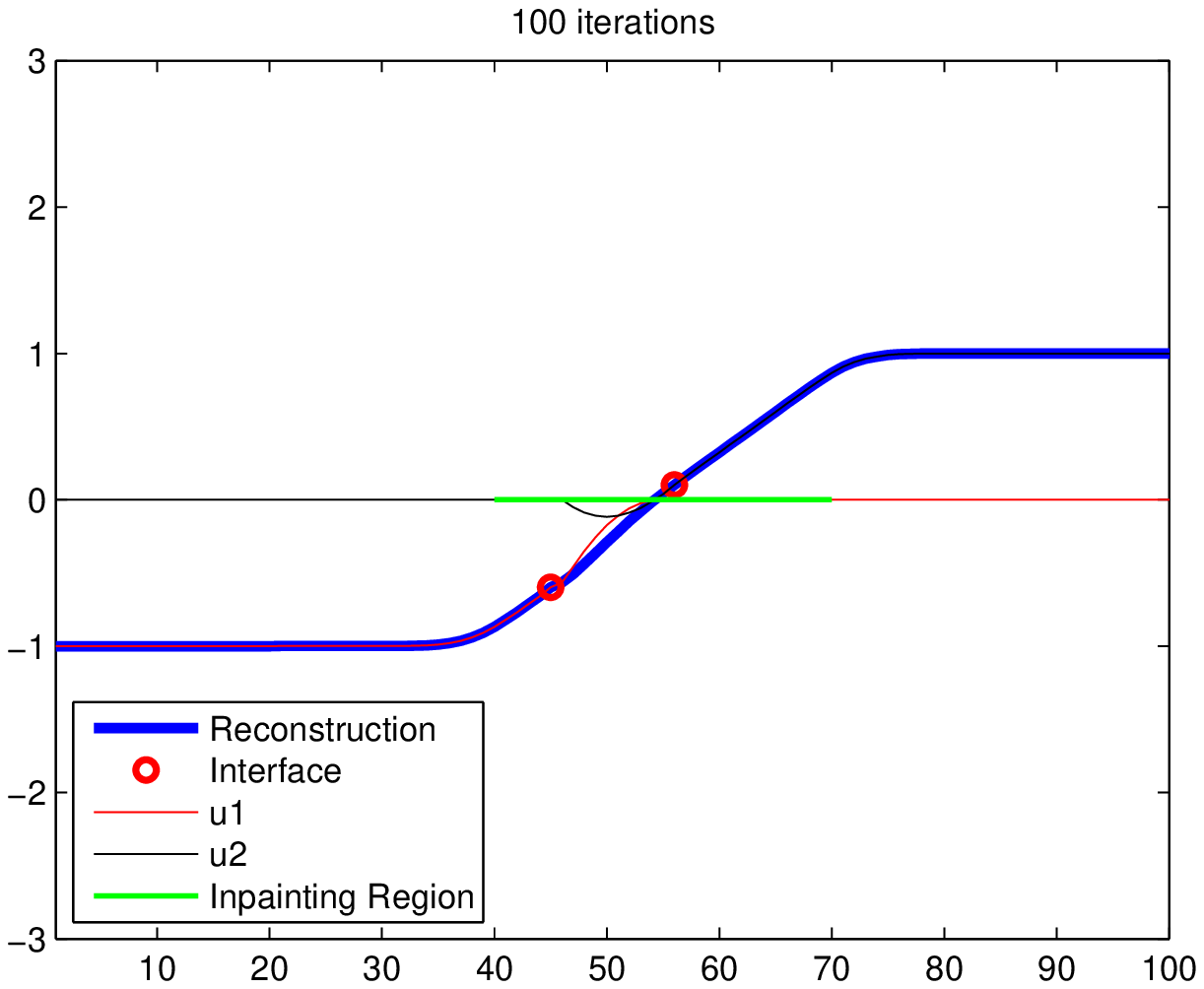}}
\end{center}    
\caption{\small Here we present two numerical experiments related to the interpolation of a 1D signal by total variation minimization, provided only information only out of an interval (indicated in green color in the figures). On the left we show an application of algorithm \eqref{schw_sp:it21} when no correction with the partition of unity is provided. In this case, the sequence of the local iterations $u_1^{(n)}, u_2^{(n)}$ is unbounded. On the right we show an application of algorithm \eqref{schw_sp:it21} with the use of the partition of unity which enforces also the uniform boundedness of the local iterations $u_1^{(n)}, u_2^{(n)}$.}
\label{fig:1Dnum}
\end{figure}

A few technical tricks are additionally required for the boundedness of the iterations in algorithm \eqref{schw_sp:it21} with respect to the nonoverlapping version. First of all the local minimizations are restricted to functions which vanish on the internal boundaries. Moreover since $u^{(n)}$ is formed as a sum of local components $u_1^{(n)}, u_2^{(n)}$ which are not uniquely determined on the overlapping part, we introduced a suitable correction by means of the partition of unity  $\{\chi_j:j=1,2\}$ in order to enforce the uniform boundedness of the sequences of the local iterations  $u_1^{(n)}, u_2^{(n)}$. With similar minor modifications, we can analogously formulate a parallel version of this algorithm as in \eqref{schw_sp:it4}.
\subsection{The solution of the local iterations}

The inner iterations 
$$
u_j^{(n+1,\ell+1)} =  \operatorname{argmin}_{\stackrel{u_j \in V_j}{\Gamma_j u_j =0}} \mathcal J_j^s(u_j+  u_{\check j}^{(n)}, u_j^{(n+1,\ell)}),
$$
where $\Gamma_j u_j =0$ is some linear constraint, are crucial for the concrete realizability of the algorithm. (Note that in the case of the nonoverlapping decomposition there is no additional linear constraint, whereas for the overlapping case we ask for the trace condition $\Gamma_j u_j = u_j|_{\Gamma_j} =0$.) Such iteration can be explicitely computed 
\begin{eqnarray*}
u_j^{(n+1,\ell+1)} &=& \left [ I- \pi_{\alpha K} \right ]\left (u_j^{(n+1,\ell)}+\pi_{V_j}T^*(g-Tu_{\check j}-T u_j^{(n+1,\ell)}) \right .\\
&& \left . \phantom{XXXXXXXXX} + u_{\check j} -\eta^{(n+1,\ell)} \right )-u_{\check j},
\end{eqnarray*}
for a suitable Lagrange multiplier $\eta^{(n+1,\ell)}$ which has the role of enforcing the linear constraints $u_j \in V_j$ and $\Gamma_j u_j =0$; $\eta^{(n+1,\ell)}$ can be approximated by an iterative algorithm, see \cite[Proposition 4.6]{fosc08} for details. Note that we have to implement repeatedly the projection $\pi_{\alpha K}$ for which the Chambolle's algorithm \eqref{chprojit} is used. More efficient algorithms can also be used such as iterative Bregman distance methods \cite{breg} or Nesterov's algorithm \cite{ne05}.

\section{Convergence properties}
 These algorithms share common convergence properties, which are listed in the following theorem.\\
{\bf Theorem.} (Convergence properties)
{\it 
The algorithms \eqref{schw_sp:it2}, \eqref{schw_sp:it4}, and \eqref{schw_sp:it21}  produce a sequence $(u^{(n)})_{n\in \mathbb{N}}$ in  $\mathcal H$ with the following properties:

(i)$\mathcal{J}(u^{(n)}) > \mathcal{J}(u^{(n+1)})$ for all $n \in \mathbb{N}$ (unless $u^{(n)}= u^{(n+1)}$);

(ii) $\lim_{n \to \infty} \| u^{(n+1)} -  u^{(n)}\|_{2} =0$;

(iii) the sequence  $(u^{(n)})_{n\in \mathbb{N}}$ has subsequences which converge in $\mathcal H$; if $(u^{(n_k)})_{k \in \mathbb{N}}$ is a converging subsequence, and  $u^{(\infty)}$ is its limit, then $u^{(\infty)}$ is always a minimizer of $\mathcal J$ in the case of algorithm \eqref{schw_sp:it21} (overlapping case), whereas for algorithms \eqref{schw_sp:it2}, \eqref{schw_sp:it4} (sequential and parallel nonoverlapping cases) this can be ensured under certain sufficient technical conditions, see \cite[Theorem 5.1 and Theorem 6.1]{fosc08} for details.  
}

\section{Numerical experiments}

\begin{figure}[ht!]
\begin{center}
\graphicspath{{./}}
    \subfigure[]{\label{l13}\includegraphics[height=3.4cm]{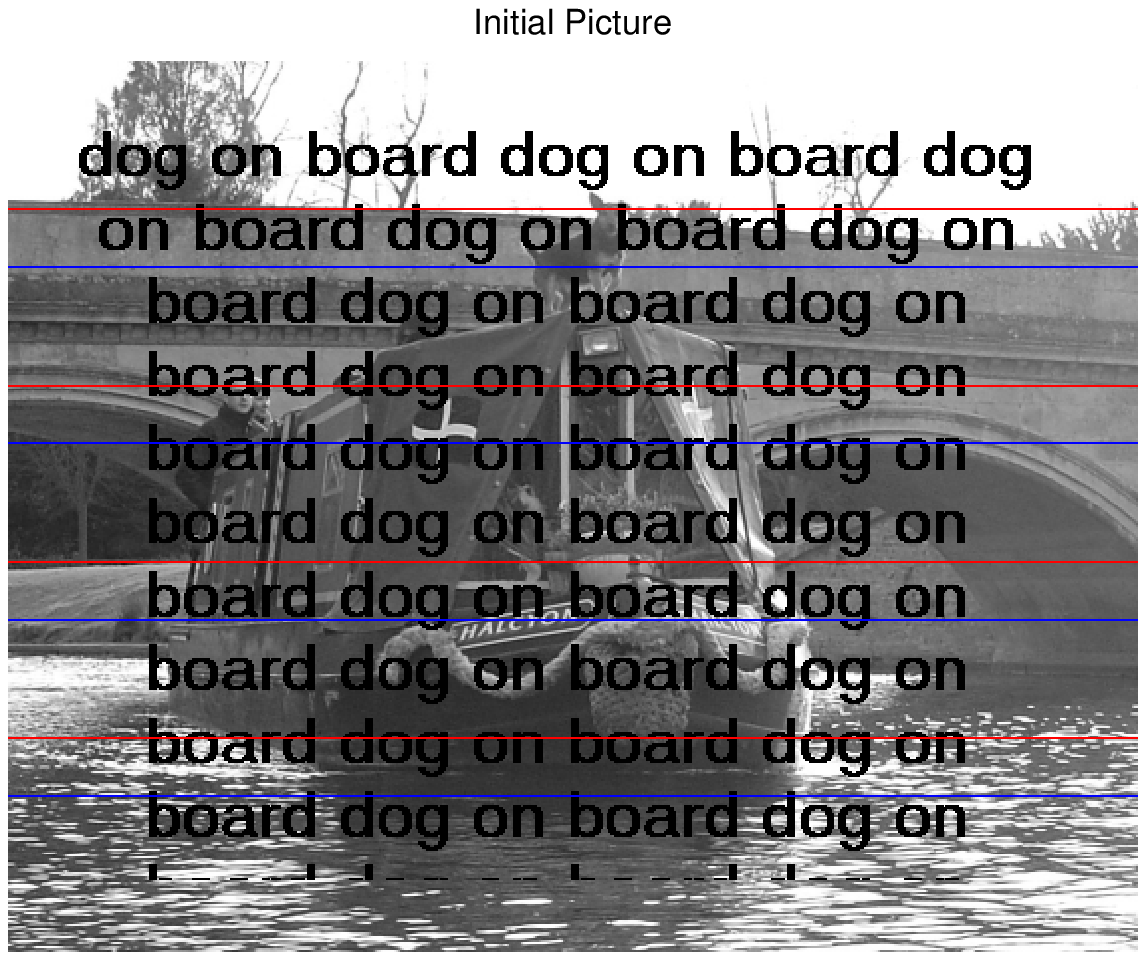}}
    \subfigure[]{\label{l23}\includegraphics[height=3.4cm]{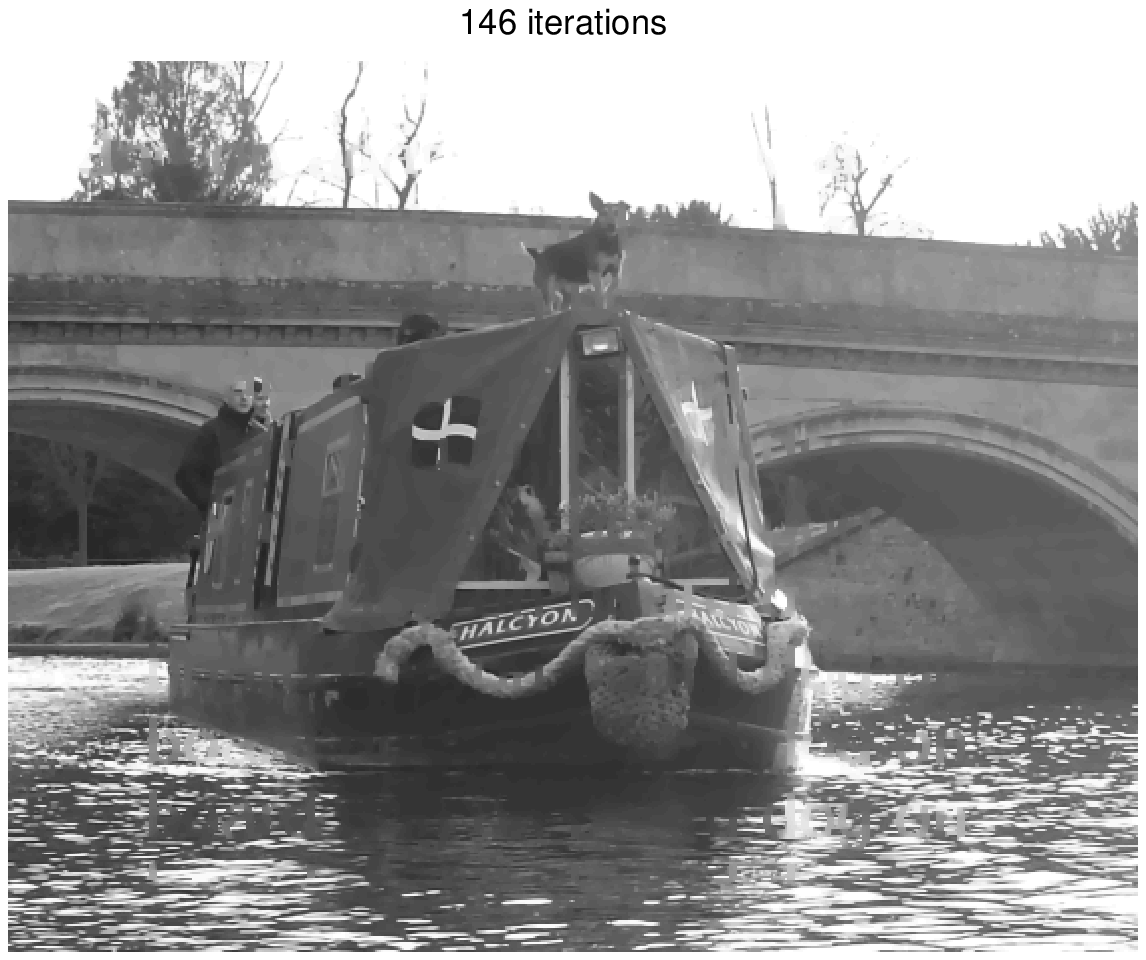}}\\
\end{center}    
\caption{\small This figure shows an application of algorithm \eqref{schw_sp:it21} for an inpainting problem. In this simulation the problem was split via decomposition into four overlapping subdomains.}
\label{fig:2Dinpainting}
\end{figure}

\begin{figure}[ht!]
\begin{center}
\graphicspath{{./}}
    \subfigure[]{\label{l11}\includegraphics[height=3.4cm]{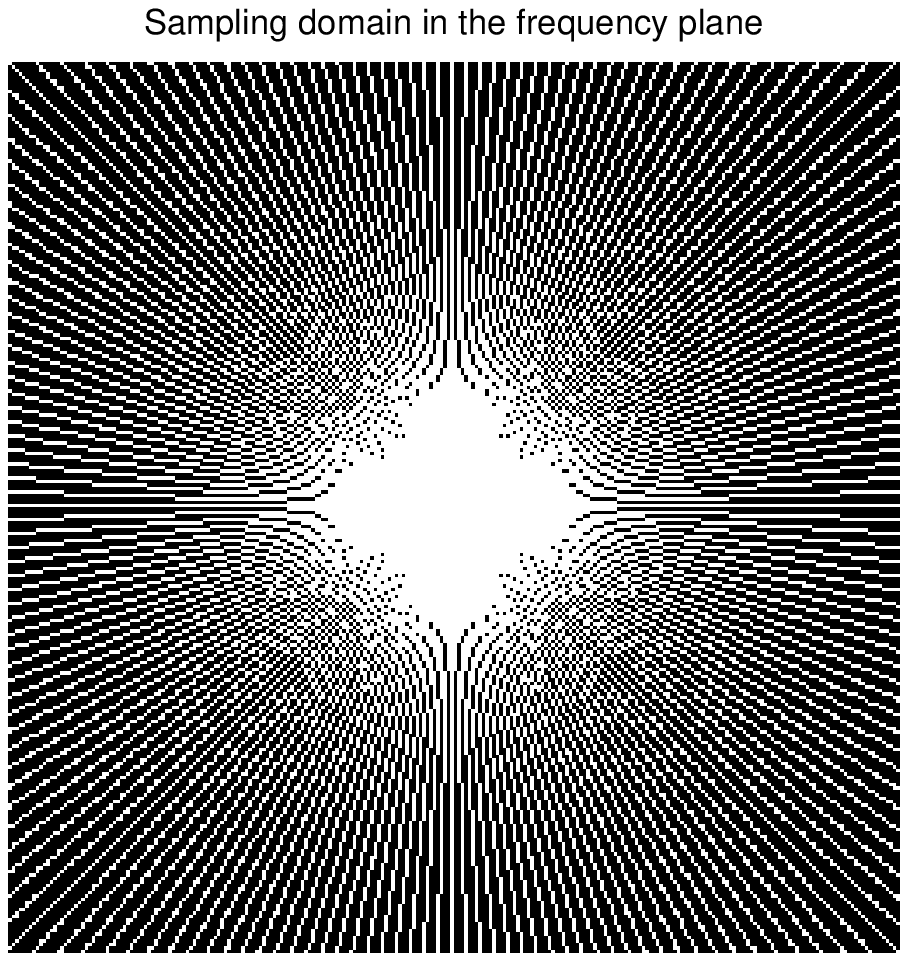}}
    \subfigure[]{\label{l21}\includegraphics[height=3.4cm]{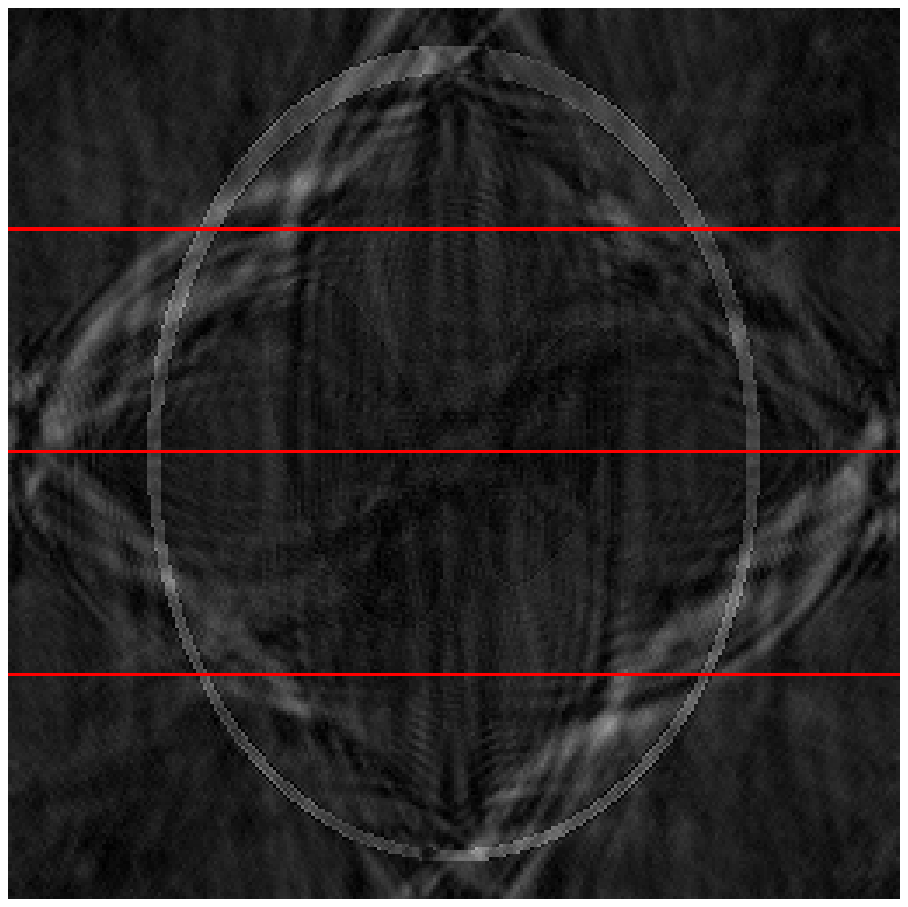}}\\
 \subfigure[]{\label{l12}\includegraphics[height=3.4cm]{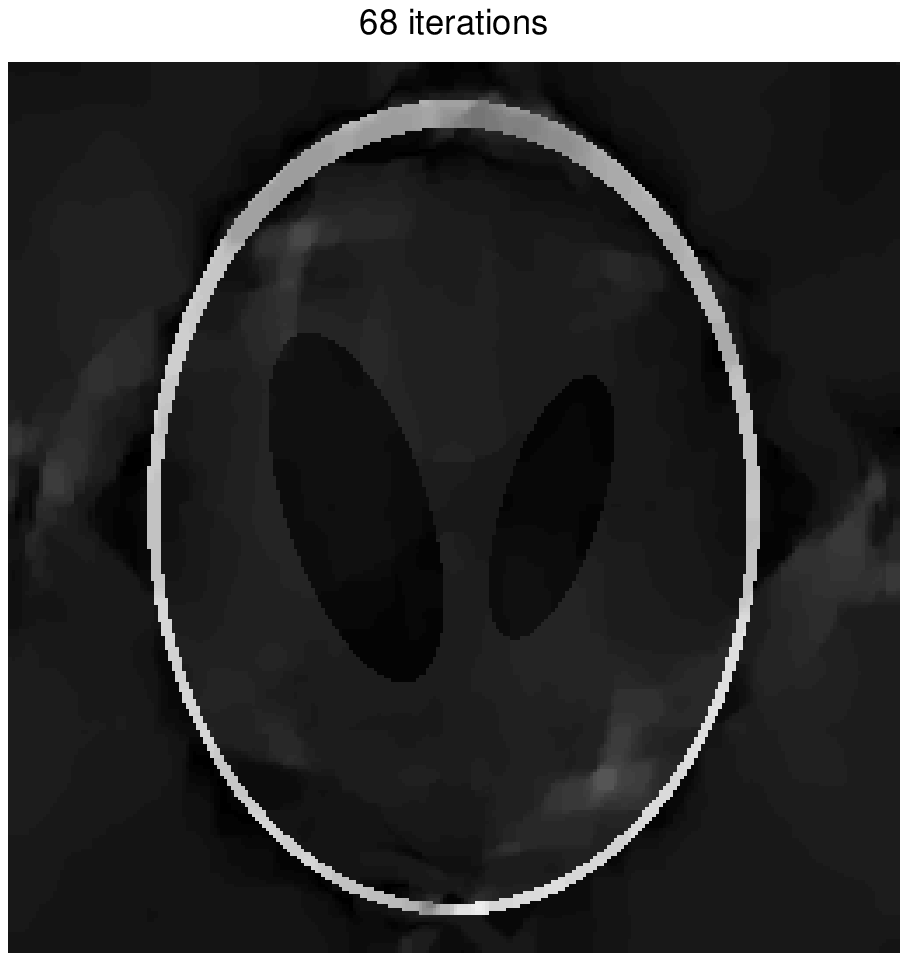}}
    \subfigure[]{\label{l22}\includegraphics[height=3.4cm]{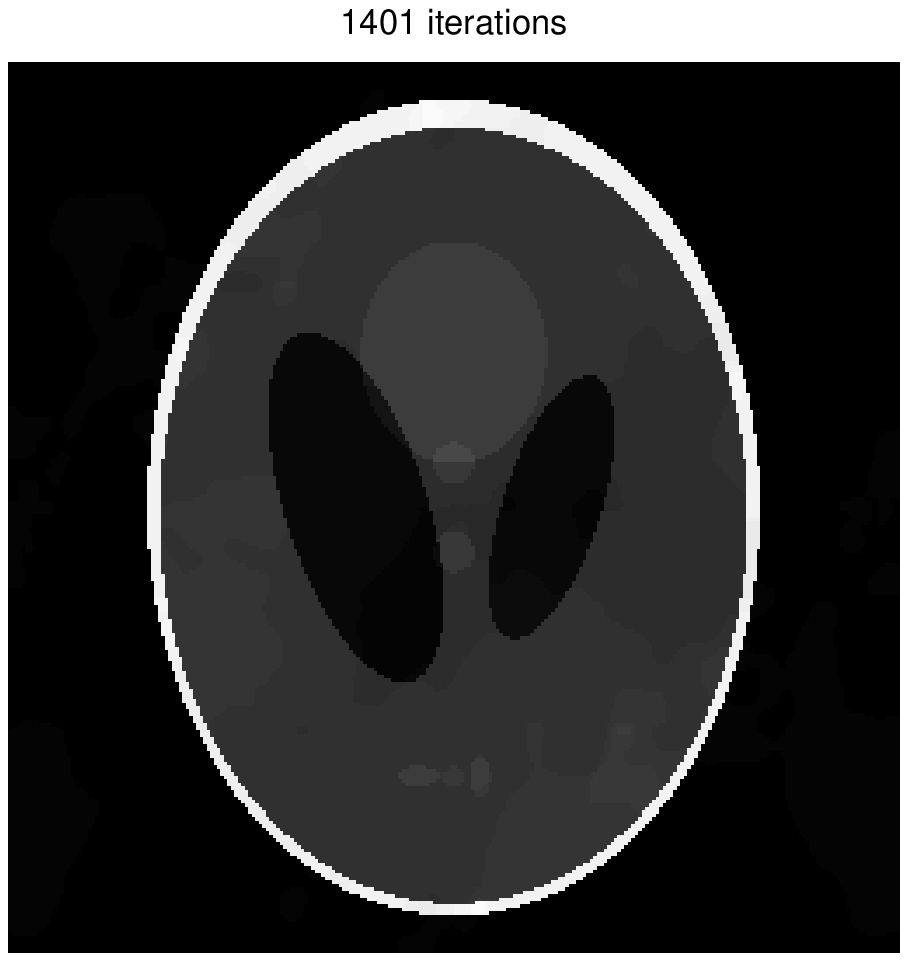}}
\end{center}    
\caption{\small We show an application of algorithm \eqref{schw_sp:it2} in a classical compressed sensing problem for recovering piecewise constant medical-type images from given partial Fourier data. In this simulation the problem was split via decomposition into four nonoverlapping subdomains. On the top-left figure, we show the sampling data of the image in the Fourier domain. On the top-right we show the back-projection provided by the sampled frequency data together with the highlighted partition of the physical domain into four subdomains. The bottom figures present intermediate iterations of the algorithm.}
\label{fig:compressedsensing}
\end{figure}

In the Figure \ref{fig:1Dnum}, Figure \ref{fig:2Dinpainting}, and Figure \ref{fig:compressedsensing} we illustrate the results of several numerical experiments, showing the successful application of algorithms \eqref{schw_sp:it2} and \eqref{schw_sp:it21}, for the restoration of 1D and 2D signals in interpolation/inpainting problems respectively, and for a compressed sensing problem.

{\footnotesize
\section*{Acknowledgments}

The authors  acknowledge the support by the project WWTF
Five senses-Call 2006, {\it Mathematical Methods for Image Analysis and
Processing in the Visual Arts}. 
M. Fornasier acknowledges the financial support provided by START grant ``Sparse Approximation and Optimization in High Dimensions'' of  the Austrian
Science Fund.
C.-B. Sch\"onlieb acknowledges the
financial support provided by KAUST (King Abdullah University of Science
and Technology), by the Wissenschaftskolleg (Graduiertenkolleg, Ph.D.
program) of the Faculty for Mathematics at the University of Vienna
(funded by the Austrian Science Fund FWF), and by the FFG project {\it Erarbeitung
neuer Algorithmen zum Image Inpainting}, projectnumber 813610.
}

{\small
\providecommand{\bysame}{\leavevmode\hbox to3em{\hrulefill}\thinspace}
\providecommand{\MR}{\relax\ifhmode\unskip\space\fi MR }
\providecommand{\MRhref}[2]{%
  \href{http://www.ams.org/mathscinet-getitem?mr=#1}{#2}
}
\providecommand{\href}[2]{#2}

\begin {spacing}{0.9}

\end{spacing}
}

\begin{thebibliography}{10}


\bibitem{AK02}
G.~Aubert and P.~Kornprobst, \emph{{M}athematical {P}roblems in {I}mage
  {P}rocessing. {P}artial {D}ifferential {E}quations and the {C}alculus of
  {V}ariation}, Springer, 2002.

\bibitem{carotaXX}
{E}.~{J}. {C}and{\'e}s, {J}. {R}omberg, and {T}. {T}ao, \emph{{E}xact signal
  reconstruction from highly incomplete frequency information}, {I}{E}{E}{E}
  {T}rans. {I}nf. {T}heory \textbf{52} (2006), no.~2, 489--509.

\bibitem{Ch} A.~Chambolle,  \emph{An algorithm for total variation minimization and applications}.  J. Math. Imaging Vision  \textbf{20}  (2004),  no. 1-2, 89--97.

\bibitem{ChL}
A.~Chambolle and P.-L. Lions, \emph{{Image recovery via total variation
  minimization and related problems.}}, Numer. Math. \textbf{76} (1997), no.~2,
  167--188.

\bibitem{CW}
P.~L. Combettes and V.~R. Wajs, \emph{Signal recovery by proximal forward-backward splitting}, Multiscale Model. Simul., \textbf{4} (2005), no.~4, 1168--1200.

\bibitem{DDFG}
I. Daubechies, R. DeVore, M. Fornasier, and S. G\"unt\"urk, {\it Iteratively re-weighted least squares minimization for sparse recovery}, Commun. Pure Appl. Math. (2009) to appear, arXiv:0807.0575

\bibitem{dateve06}
I.~Daubechies, G.~Teschke, and L.~Vese, \emph{Iteratively solving linear
  inverse problems under general convex constraints}, Inverse Probl. Imaging
  \textbf{1} (2007), no.~1, 29--46.

\bibitem{do04}
{D}.~{L}. {D}onoho, \emph{{C}ompressed sensing}, {I}{E}{E}{E} {T}rans. {I}nf.
  {T}heory \textbf{52} (2006), no.~4, 1289--1306.

\bibitem{fo07}
 M.~Fornasier, \emph{Domain decomposition methods for linear inverse problems with sparsity constraints},
  Inverse Problems \textbf{23} (2007), no. 6, 2505--2526.

\bibitem{fosc08} 
M.~Fornasier and C.-B. Sch\"onlieb, {\it Subspace correction methods for total variation and $\ell_1-$minimization}, arXiv:0712.2258 (2007)


\bibitem{ne05} Y.~Nesterov, {\it Smooth minimization of non-smooth functions}. Mathematic Programming,
Ser. A, \textbf{103} (2005), 127--152.

\bibitem{breg}
S. Osher, M. Burger, D. Goldfarb, J. Xu, and W. Yin, {\it An Iterative Regularization Method for Total Variation-Based Image Restoration},
Multiscale Model. Simul. \textbf{4}, no.~2 (2005) 460-489.

\bibitem{ROF}
L.~I. Rudin, S.~Osher, and E.~Fatemi, \emph{{Nonlinear total variation based
  noise removal algorithms.}}, Physica D \textbf{60} (1992), no.~1-4, 259--268.

\bibitem{Ve01}
L.~Vese, \emph{{A study in the BV space of a denoising-deblurring variational
  problem.}}, Appl. Math. Optim. \textbf{44} (2001), 131--161.

\bibitem{WBA}
P. Weiss, L.~Blanc-F\'eraud, and G. Aubert, {\it Efficient schemes for total variation minimization under constraints in image processing}, SIAM J. Sci. Comput., (2009) to appear.
\end{thebibliography}
\end{document}